\documentclass[12pt]{article}

\usepackage{amssymb}
\usepackage[leqno]{amsmath}
\usepackage{amsthm}

\newtheorem{theorem}{Theorem}[section]

\newtheorem{lemma}[theorem]{Lemma}

\theoremstyle{definition}

\newcommand{\RR}{{\mathbb R}}

\begin{document}

\title{Isotopies of high genus Lagrangian surfaces}

\author{R. Hind \thanks{Supported in part by NSF grant DMS-0505778.} \and A. Ivrii}

\date{\today}

\maketitle

\section{Introduction}

In this paper we address a local version of the isotopy problem
for Lagrangian surfaces in a symplectic $4$-manifold $(M,\omega)$.
This question was first raised by V. Arnold in \cite{arnold}. A
Lagrangian submanifold $L$ is one for which $\omega|_L$ vanishes.
In general we would like to classify homotopic Lagrangian
submanifolds up to smooth isotopy or better still Lagrangian
isotopy, that is, smooth isotopy through Lagrangian submanifolds.
Equivalence classes are called Lagrangian knots. Here we show that
in a sufficiently small neighborhood of a given Lagrangian surface
there are no Lagrangian knots up to smooth isotopy. More precisely
our result can be stated as follows.

\begin{theorem}
Let $T^* \Sigma$ be the cotangent bundle of a Riemann surface with
its canonical symplectic structure and $L \subset T^* \Sigma$ be a
connected Lagrangian submanifold homologous to $\Sigma$. Then $L$
is smoothly isotopic to $\Sigma$.
\end{theorem}

In the case when $\Sigma$ has genus $0$ or $1$ the above Theorem
1.1 is due to Y. Eliashberg and L. Polterovich, see \cite{elpolt}.
In fact, work of the first author, see \cite{hind}, shows that if
$\Sigma$ has genus $0$ then all such Lagrangian spheres are
Lagrangian isotopic to the zero-section. However we remark that it
is not true in general that isotopic spheres are Lagrangian
isotopic, see the work of P. Siedel, \cite{seidel}. If $\Sigma$
has genus $1$ work of the second author \cite{ivrii} shows again
that all such $L$ are Lagrangian isotopic. The question of whether
or not in higher genus cases all such isotopic Lagrangians are
Lagrangian isotopic remains open.

In general symplectic $4$-manifolds there exist homologous high
genus Lagrangian submanifolds which are not smoothly isotopic, see
the work of D. Park, M. Poddar and S. Vidussi, \cite{park}.

\section{Proof of the theorem}

In this section we prove Theorem 1.1. Since the result is known
when $\Sigma$ has genus $0$ or $1$ we will assume throughout that
$\Sigma$ has genus $g
>1$.
Let $\sigma$ be an area form on $\Sigma$ of total area $2g-2$. Let
$\pi:T^* \Sigma \to \Sigma$ be the projection along the fibers.
The cotangent bundle $T^* \Sigma$ carries a canonical symplectic
form $\omega_0=d(\lambda_0)$, where $\lambda_0=pd \pi$ is the
Liouville form. The zero section $\Sigma$ is Lagrangian with
respect to $\omega_0$.

We can also think of $T^* \Sigma$ as a tubular neighborhood of a
symplectic submanifold $\Sigma$. Then $T^* \Sigma$ carries another
symplectic form $\tau$ which is symplectic on the fibers and such
that $\tau|_\Sigma =\sigma$. Let $r:T^* \Sigma \to [0,\infty)$ be
the length function with respect to an Hermitian metric on (the
complex line bundle) $T^* \Sigma$. We denote the levels by $T^r
\Sigma$. Then the unit circle bundle $\pi :T^1 \Sigma \to \Sigma$
carries a connection $\alpha$ with $d\alpha = \pi^* \sigma$. We
can arrange that $\tau|_{T^r \Sigma} = f(r)d\tilde{\alpha}$ where
$\tilde{\alpha}$ is the pullback of the form $\alpha$ on $T^1
\Sigma$ and $f$ is decreasing towards $0$ as $r$ approaches
$\infty$.

For $\epsilon$ sufficiently small,
$\Omega_{\epsilon}=\omega_0+\epsilon \tau$ is also a symplectic
form on $T^* \Sigma$.

We reparameterize $\omega_0$ such that outside of a large compact
set $T^{\le r_0} \Sigma$ it is given by $d(e^r \tilde{\lambda_0})$
where $\tilde{\lambda_0}$ denotes the pullback of the Liouville
form from the unit tangent bundle. Also outside of $T^{\le r_0}
\Sigma$ we extend $\tau$ by extending the function $f$ to a
decreasing function $g(r)$ with $g=-e^r$ outside of a (larger)
compact set. Then we define a new form $\omega$ on $T^* \Sigma$ by
$\omega = \Omega_{\epsilon}$ on $T^{\le r_0} \Sigma$ and $\omega =
d(e^r \tilde{\lambda_0} + \epsilon g(r) \tilde{\alpha})$
elsewhere. We note that $\omega$ is a symplectic form for
$\epsilon$ sufficiently small and that the fibers of $T^* \Sigma$
are $\omega$-symplectic planes of infinte area.

Let $V$ be a tubular neighborhood of our Lagrangian submanifold $L
 \subset (T^{\le r_0} \Sigma, \omega_0)$.

\begin{lemma}
There exists an $\epsilon_0>0$ such that for all $\epsilon <
\epsilon_0$ the Lagrangian $L$ can be isotoped to an
$\Omega_\epsilon$ symplectic surface within $V$.
\end{lemma}

{\bf Proof} This is a slight modification of Proposition 2.1.A in
\cite{elpolt}. Let $\sigma$ be a symplectic form on $V$ such that
$\sigma|_L$ is an area form of total area $2g-2$. Then $(\tau -
\sigma)|_L$ is exact and so by the relative Poincar\'{e} Lemma
there exists a $1$-form $\lambda$ on $V$ such that $\sigma=\tau
+d\lambda$. Let $\rho:V \to [0,1]$ have compact support and equal
$1$ close to $L$. Then there exists an $\epsilon_0$ such that for
$\epsilon \le \epsilon_0$ the form
$\Omega_{\epsilon}'=\omega_0+\epsilon(\tau +d(\rho \lambda))$ is
symplectic as is the linear family of forms connecting
$\Omega_{\epsilon}'$ to $\Omega_{\epsilon}$. As
$\Omega_{\epsilon}'=\Omega_{\epsilon}$ away from $V$ and $L$ is
$\Omega_{\epsilon}'$ symplectic it follows from Moser's method
that $L$ can be isotoped to an $\Omega_{\epsilon}$ symplectic
surface inside $V$. \qed

The results from \cite{hindivrii} imply that that all
$\Omega_{\epsilon}$ symplectic surfaces sufficiently close to
$\Sigma$ are isotopic to $\Sigma$. Thus we could conclude here if
it were possible to arrange that the symplectic surface was
contained in a suitably small symplectic neighborhood. However we
could find no straightforward method of doing this. Instead we
proceed as follows.

\begin{lemma} For $\epsilon$ sufficiently small, all
connected symplectic surfaces $S$ in $(T^* \Sigma, \omega)$ which
are homologous to $\Sigma$ and intersect the fiber over a point
$p$ exactly once transversally must be smoothly isotopic to
$\Sigma$.
\end{lemma}

{\bf Proof} Let $U$ be a neighborhood of $p$ such that a given
symplectic surface $S$ intersects all fibers over points $q \in U$
transversally in a single point. By a small perturbation we may
assume that $S \cap \Sigma$ is disjoint from $\pi^{-1}(U)$.

Let $h:\Sigma \to \RR$ be a Morse function with a single minimum
and all critical points contained in $U$. Then the gradient
flowlines of $h$ foliate the complement of the critical points of
$h$ by curves $\gamma(x):(-\infty, \infty) \to \Sigma$ which lie
in $U$ for $|x|$ sufficiently large. Denote the critical points of
$h$ by $p_1,...,p_N$.

Let $s_i \in S$ be the unique point with $\pi (s_i)=p_i$. Then we
also assume that as subspaces of $T(T^* \Sigma)$ we have $T_{s_i}
S = T(\pi^{-1}(p_i))^{\perp \Omega_{\epsilon}}$, the symplectic
complement to the tangent space of the fiber.

Recall that for $r$ sufficiently large $\omega|_{T^r \Sigma} =
d\beta$ where $\beta = e^r \tilde{\lambda_0} + \epsilon g(r)
\tilde{\alpha}$ is a contact form for $\epsilon$ sufficiently
small. We observe that $\pi^{-1}(\gamma(\RR))\cap T^r \Sigma$ is a
cylinder $C_{\gamma}$ foliated by the circles $F_x=
\pi^{-1}(\gamma(x))$. Now, $\tilde{\lambda_0}$ vanishes on the
$F_x$ while $\tilde{\alpha}$ does not. Therefore $\ker
\beta|_{C_{\gamma}}$ is a nonsingular line field $l$ transverse to
all $F_x$. In particular $l$ has no closed orbits.

We claim that there exists an almost complex structure $J_0$ on
$T^* \Sigma$ which is tamed by $\omega$ and satisfies the
following properties. The surfaces $S \cap \pi^{-1}(U)$ and
$\Sigma$ are $J_0$ holomorphic; the contact planes $\ker \beta$ in
$T^r \Sigma$ are $J_0$ holomorphic for some $r$ sufficiently
large; for all critical points $p_i$ the disk $D_i =
\pi^{-1}(p_i)$ is $J_0$-holomorphic.

The only requirement here which is not well known is the claim
that it is possible to find a $J_0$ along $T^r \Sigma$ which
simultaneously makes both the subbundles $\ker \beta$ and
$\pi^{-1}(p_i)$ into $J_0$-holomorphic distributions. But the
existence of such $J_0$ is established in a more general context
by Theorem 7.4 in the article \cite{coffey} of J. Coffey.

Let $J_t$, $0 \le t \le 1$ be a family of almost-complex
structures on $T^* \Sigma$ coinciding with $J_0$ outside some $T^s
\Sigma$, where $s<r$, and on $\pi^{-1}(U)$, such that $S$ is $J_1$
holomorphic.

We next claim that for all $t$ the cylinders $C_{\gamma}$ can be
foliated by circles which bound $J_t$ holomorphic disks. These
circles are transverse to $l$ and at the ends of the cylinders the
holomorphic disks converge to the perturbed fibers $\pi^{-1}(p_i)$
for $p_i$ a critical point. The union of all disks over all
cylinders gives a foliation of $T^{\le r} \Sigma$ by disks in the
relative homotopy class of the fibers.

This claim follows from the theory of filling by holomorphic
disks, see \cite{filling}. For each $\gamma$, the cylinder
$C_{\gamma}$ is foliated by the boundaries of embedded holomorphic
disks near its ends. But as the cylinders are totally real the
foliation extends to cover the whole cylinder. The only
obstruction in this case is bubbling of holomorphic spheres inside
$T^* \Sigma$ and bubbling of disks on the boundary. But as
$\pi_2(\Sigma)$ is trivial such spheres do not exist. Bubbling of
disks can be excluded as in \cite{filling} since all holomorphic
disks with boundary on $T^r \Sigma$ must have boundary transverse
to $l$. For embedded boundaries this fixes the homology class and
prevents degeneracies.

The disks $D_i$ constructed above are $J_t$ holomorphic for all
$t$ and their intersection with $S$ and $\Sigma$ is transversal
and in a single point. Therefore by positivity of intersections
the same is true for all intersections of $J_0$ holomorphic disks
with $\Sigma$ and all $J_1$ holomorphic disks with $S$.

We fix a Riemannian metric on $T^* \Sigma$ which decays rapidly
along the fibers. Then with respect to the restricted metric the
centers of mass of our holomorphic disks give a smooth family of
surfaces $G_t$. By the previous remark, it is clear that $G_0$ is
smoothly isotopic to $\Sigma$ and  $G_1$ is isotopic to $S$. \qed

\begin{lemma} The Lagrangian $L$ can be isotoped to an $\omega$ symplectic surface in $T^* \Sigma$ intersecting the fiber over a point $p$ transversally in a single point.
\end{lemma}

This will be established using the theory of finite energy planes.
This theory was developed in the series of papers \cite{hof},
\cite{hofa}, \cite{hofi}, \cite{hoff}. Based on these, work of the
first author \cite{hind1}, \cite{hind2} deduced the existence of
finite energy planes in $T^* S^2$ lying in certain relative
homotopy classes. The result we need is the following.

\begin{theorem}
There exists a Morse-Bott type contact form on the unit cotangent
bundle $T^1 S^2$ with an isolated Reeb orbit $\gamma$ having
minimal action and Conley-Zehnder index $1$. The form can be
chosen arbitrarily close to the standard form (where all orbits
have action $2\pi$).

The orbit $\gamma$ is the asymptotic limit of exactly two finite
energy planes in $T^* S^2$ which have intersection number $\pm 1$
with the zero-section.
\end{theorem}

Here we think of $T^* S^2$ as an almost-complex symplectic
manifold with convex end corresponding to the positive
symplectization of $T^1 S^2$.

{\bf Outline of proof of Theorem $2.4$}

Our contact form will be a perturbation of the standard contact
form on $T^1 S^2$ defined with the round metric. With respect to
the standard form the space of closed Reeb orbits can be
identified with a $2$-sphere. We now perturb the contact form
following \cite{B} such that the resulting form is of Morse-Bott
type. We do this using a function on the orbit space of closed
Reeb orbits. We choose a function which is rotationally symmetric
and has critical points at two isolated critical points and at a
finite number of radii. The critical points on the orbit space
correspond to closed Reeb orbits for the perturbed contact form,
so we obtain two isolated Reeb orbits and a number of
$1$-parameter families. Let $\gamma$ be an isolated Reeb orbit
corresponding to a miniumum of our function which therefore is of
minimal period as required.

We remark here that equivalently our contact form can be realized
as the restriction of the Liouville form on $T^* S^2$ to a
hypersurface which is a suitable perturbation of the unit
cotangent bundle for the round metric.

To produce finite energy planes asymptotic to the orbit of minimal
length we follow the analysis in \cite{hind1}, see also
\cite{hind2}. Our contact manifold $T^1 S^2$ is doubly covered by
the contact $S^3$. Each of the closed orbits is doubly covered by
a closed orbit of the corresponding flow on $S^3$. Now the method
of filling by holomorphic disks can be applied as in \cite{hof} to
produce a finite energy plane in the symplectization $\Bbb R
\times S^3$ asymptotic to a given periodic orbit of minimal
length. Such a finite energy plane projects to give a finite
energy plane in $T^* S^2$ asymptotic to twice the corresponding
simple orbit.

In $T^* S^2$ the (unparameterized) finite energy planes asymptotic
to twice a simple orbit appear in a $2$-dimensional family which
foliates a region of $T^* S^2$. Now the considerations in
\cite{hind1} imply that this family cannot be compact and a
certain subsequence of such planes will bubble to produce a pair
of planes asymptotic to the simple orbit. The planes have
intersection $\pm 1$ with the zero-section as required.

There exists a circle of isometries of the round metric on $S^2$
which preserve the geodesic corresponding to $\gamma$ (simply the
rotations about a perpendicular axis). The isometries lift to a
family of symplectomorphisms of $T^* S^2$ restricting to a family
of contactomorphisms of $T^1 S^2$ preserving the Reeb orbit
$\gamma$. We may assume that this family of symplectomorphisms
also preserves the perturbed $T^1 S^2$ and the closed orbit
$\gamma$. If the almost-complex structure we use is also invariant
under this family of rotations then our planes asymptotic to
$\gamma$ must also be invariant. Their uniqueness follows as in
\cite{hinds}.

 \qed

 {\bf Remark}
 At least if the family is also of minimal length then the same
 construction can be used to produce a finite energy plane
 asymptotic to a closed Reeb orbit in one of the $1$-parameter
 families surrounding $\gamma$ in $T^1 S^2$. Acting by the circle
 of isometries we obtain finite energy planes asymptotic to each
 orbit in this family. The finite energy plane asymptotic to
 $\gamma$ has deformation index $0$. Further, arguing similarly to
 Theorem 4.3 in \cite{hofi}, it can be shown that the
 corresponding Cauchy-Riemann operator is surjective. It seems to
 be a subtle question to decide whether the finite energy planes
 in the $1$-parameter family are also regular, but we do not need
 this fact in what follows.

We can think of choosing a contact form on $T^1 S^2$ as being
equivalent to choosing a Finsler metric on $S^2$, the Reeb orbits
then correspond to Finsler geodesics. It is also equivalent to
choosing a length function in $T^* \Sigma$ and restricting the
Liouville form to the level sets.

Now, in the above construction, the projections to the
zero-section of the two planes asymptotic to $\gamma$ can be
compactified to maps from a closed disk with boundary on the
corresponding (Finsler) geodesic. For almost-complex structures
sufficiently close to the standard one the images of these disks
should occupy opposite hemispheres in the $S^2$. In any case, if
we restrict attention to one relative homotopy class, say that
with planes having intersection $+1$ with the zero-section, then
we may assume that the projection to $S^2$ of the plane asymptotic
our minimal orbit is disjoint from a segment $\sigma_0$ of a
geodesic $\sigma$ which intersects the geodesic corresponding to
$\gamma$ orthogonally at two points.

We may now assume that there exists an $\epsilon>0$ such that
$\gamma$ and the geodesics in the surrounding $1$-parameter family
have length less than $2\pi - \epsilon$ and all other closed
geodesics have length at least $2\pi + \epsilon$. In particular
$\sigma_0$ can be taken to have length greater than $\pi$. We can
further arrange for $\sigma_0$ to be disjoint from the
$1$-parameter family.

{\bf Proof of Lemma $2.3$} We choose a Finsler metric on $\Sigma$
and fix an embedded geodesic segment $\sigma_1$. Scaling the
metric appropriately we can increase the length of this segment to
$\pi$ and assume that any closed geodesics on $\Sigma$ or geodesic
segments which start and end on $\sigma_1$ have length at least
$2\pi$.

On our $S^2$ from Theorem $2.4$ we choose neighborhoods $U$ of
$\sigma_0$ and $V$ of $\sigma$ such that the following conditions
hold. Geodesic segments starting and ending in $U$ either lie
entirely in $V$ or have length at least $2\pi$; closed
noncontractible curves in $V$ have length at least $2\pi$; closed
noncontractible curves in $V \setminus \sigma_0$ have length at
least $2\pi$.

We now perform a connected sum of the Finsler surfaces $\Sigma$
and the $S^2$ to obtain a new Finsler metric on a surface
diffeomorphic to $\Sigma$ as follows. We do this by removing very
small neighborhoods $D_0$ and $D_1$ of $\sigma_0$ and $\sigma_1$
respectively and replacing these with a cylinder. We can extend
the Finsler metric over the cylinder so that curves traversing its
longitude must have length at least $6\pi$ and noncontractible
loops in the cylinder have length at least $2\pi$. Furthermore, we
may assume that the union of $S^2 \setminus D_0$ and one half of
the cylinder glued onto it can be isometrically embedded in $S^2
\setminus \sigma_0$ with a new complete Finsler metric, the
embedding being the identity on $S^2 \setminus D_0$. This complete
Finsler metric evaluated on tangent vectors is everywhere bounded
below by the original metric on the $S^2$.

We claim that $\gamma$ is the closed geodesic of minimal length on
our new surface. This follows from the construction for geodesics
disjoint from the glued cylinder and for geodesics intersecting
both the original $S^2$ and the original $\Sigma$. It is easily
guaranteed that the cylinder has a longitudinal foliation by
geodesic segments and so geodesics lying in entirely in the
cylinder must be noncontractible and so also have length at least
$2\pi$. Suppose then that a geodesic of length less than $2\pi$
lies partly in the cylinder and partly, say, in the original
$S^2$. Then it can be identified with a geodesic in $S^2 \setminus
\sigma_0$ with its complete Finsler metric. Thus if it stays in
the region $V$ it must again have length at least $2\pi$, but if
it leaves $V$ then it contains a geodesic segment starting and
ending in $U$ which also must have length at least $2\pi$. This
justifies our claim.

We now vary the symplectic form and tame almost-complex structure
inside a compact set so that eventually the form is equal to
$\Omega_{\epsilon}$ near the zero-section and a symplectic surface
$S$ isotopic to $L$ is holomorphic.

Now, for a suitable almost-complex structure $J_0$ the finite
energy plane from Theorem $2.4$ asymptotic to $\gamma$ and having
intersection number $+1$ with the $S^2$ also exists in the
cotangent bundle of the new surface $\Sigma$. Indeed, it lies
disjoint from the region where the connected sum was performed.
Then by positivity of intersection with the $1$-parameter family
of planes asymptotic to the orbits surrounding $\gamma$ we see
that any finite energy plane asymptotic to $\gamma$ in this
homotopy class must lie in the part of $T^* \Sigma$ projecting to
the $S^2$ and so from Theorem $2.4$ the finite energy plane is
unique.

We claim that this finite energy plane asymptotic to our Reeb
orbit $\gamma$ continues to exist as the complex structure is
deformed to a new structure $J_1$. By positivity of intersections
we will then see that eventually the plane must intersect $S$
transversally in a single point.

But since $\gamma$ has minimal period the space of planes
asymptotic to $\gamma$ is compact modulo reparameterizations, see
\cite{BEHWZ}. Combining this with the Fredholm theory for finite
energy planes, the space of planes holomorphic with respect to
$J_1$ is cobordant to the space of planes holomorphic with respect
to $J_0$. There is a unique (regular) $J_0$-holomorphic plane and
therefore there is an odd number of $J_1$-holomorphic planes
asymptotic to $\gamma$, which justifies our claim.

Hence our finite energy planes persist as the almost-complex
structure is deformed and we find an isotopy from the initial
plane to a plane which intersects $S$ in a single point. If
$\epsilon$ is chosen sufficiently small then the planes can both
be assumed to be symplectic with respect to $\omega$ and the
isotopy through $\omega$-symplectic planes. The forms $\omega$ on
$T^* \Sigma$ were described at the start of this section, we
recall that restricted to a level $T^r \Sigma$ we have
$\omega=d\beta$ where $\beta$ is a small perturbation of the
Liouville form.

We recall the initial plane may be assumed to intersect the
zero-section transversally in a single point and actually each
level $T^r \Sigma$ in a curve transverse to the contact structure
(see for instance \cite{hinds}). Using this we can construct an
$\omega$-symplectic isotopy taking the initial plane into a fiber
through planes similarly intersecting the levels in transverse
curves, as follows.

First we remark that given a family $\{\gamma_t\} \subset (T^{r_0}
\Sigma, \beta)$ of transverse curves for $t \ge 0$ the union
$\{(t,\gamma_t)|t \ge 0\}$ is a sympectic surface in $((0,\infty)
\times T^{r_0} \Sigma, d(e^t \beta))$ provided $\gamma_t$ changes
sufficiently slowly with $t$. (Its tangent space is spanned at a
point $p$ by the tangent to $\gamma_t$ in $\{t\} \times T^{r_0}
\Sigma$ and $\frac{\partial}{\partial t} + \frac{\partial
\gamma_t}{\partial t}(p)$.) Moreover if $\{\gamma_t\}$ already
generates a symplectic surface and we reparameterize the family
$\{\gamma_t\}$ to reduce this rate of change then the surface
remains symplectic.

To construct our isotopy, given these remarks, it suffices to find
an isotopy of transverse curves in some $(T^r \Sigma, \beta)$ from
the Reeb orbit corresponding to $\gamma$ to a fiber of $\pi$. Such
an isotopy can be constructed as follows. We note that the natural
lift (via the metric dual of its derivative) of any embedded curve
in $\Sigma$ to $T^r \Sigma$ is a transverse curve with respect to
the (Liouville) contact structure induced from the metric on
$\Sigma$. The Reeb orbit corresponding to $\gamma$ is the natural
lift of the meridian and as the meridian contracts to a point we
get an isotopy of transverse curves to a curve $C^1$ close to a
fiber. This is also a transverse isotopy with respect to $\beta$
for $\epsilon$ sufficiently small. A transverse curve sufficiently
close to a fiber can be moved into a fiber through curves
transverse with respect to $\beta$. We then construct a symplectic
isotopy from the isotopy of transverse curves by observing again
that a union of transverse curves $\gamma_s$ in the contact
manifolds $T^s \Sigma$ is a symplectic surface provided that
$\gamma_s$ changes sufficiently slowly with $s$. With respect to
$\omega$ the fibers of $T^* \Sigma$ have infinite area and so this
isotopy can be arranged to leave a neighborhood of the
zero-section inside a given compact set.

Putting everything together, we obtain a proper isotopy from a
fiber to a plane intersecting $S$ transversally in a single point.
The preimage of $S$ is a surface intersecting a fiber in a single
point as required. \qed

\vspace{0.1in}

Richard Hind\\
Department of Mathematics\\
University of Notre Dame\\
Notre Dame, IN 46556\\
email: hind.1@nd.edu

\vspace{0.1in}

Alexander Ivrii\\
D\'epartment de Math\'ematiques et de Statistique\\
Universit\'e de Montr\'eal\\
CP 6128 succ Centre-Ville\\
Montr\'eal, QC H3C 3J7, Canada\\
email: ivrii@DMS.UMontreal.CA


\newpage

\begin{thebibliography} {99}

\bibitem{arnold} V. I. Arnold, First steps in symplectic topology,
{\it Russ. Math. Surv.}, 41 (6) (1986), 1-21.
\bibitem{B}
F.~Bourgeois.
\newblock A Morse-Bott approach to contact homology. Symplectic and contact topology: interactions and perspectives (Toronto, ON/Montreal, QC, 2001), 55--77, Fields Inst. Commun., 35, Amer. Math. Soc., Providence, RI, 2003.
\newblock 2002.
\bibitem{BEHWZ}
F.~Bourgeois, Y.~Eliashberg, H.~Hofer, K.~Wysocki, and E.~Zehnder.
\newblock Compactness results in symplectic field theory.
\newblock {\em Geom. Topol.}, 7:799--888 (electronic), 2003.
\bibitem{coffey} J. Coffey, Symplectomorphism groups and isotropic
skeletons, {\it Geom. Topol.}, 9 (2005), 935-970.
\bibitem{filling} Y. Eliashberg, Filling by holomorphic discs and its applications. Geometry of low-dimensional manifolds, 2 (Durham, 1989), 45--67, London Math. Soc. Lecture Note Ser., 151, Cambridge Univ. Press, Cambridge, 1990.
\bibitem{elpolt} Y. Eliashberg and L. Polterovich,
Unknottedness of Lagrangian surfaces in symplectic $4$-manifolds,
{\it Internat. Math. Res. Notices}, 11 (1993), 295--301.
\bibitem{EGH} Y. Eliashberg, A. Givental and H. Hofer, Introduction to symplectic field theory, GAFA 2000 (Tel Aviv, 1999), {\it Geom. Funct. Anal.}, 2000, Special Volume, Part II, 560--673.

\bibitem{fint} R. Fintushel and  R. J. Stern, Symplectic surfaces in a fixed homology class, {\it J. Differential Geom.}, 52 (1999), no. 2, 203--222.
\bibitem{Gr} M. Gromov, Pseudo-holomorphic curves in symplectic manifolds, {\it Inv. Math.}, 82(1985), 307-347.
\bibitem{hind1} R. Hind, Holomorphic Filling of $\Bbb RP^3$, {\it Comm. in Contemp. Math.}, 2(2000), 349-363.
\bibitem{hind2} R. Hind, Stein fillings of Lens spaces, {\it Comm. in Contemp. Math.}, 6(2003), 1-17.
\bibitem{hinds} R. Hind, Lagrangian spheres in $S^2 \times S^2$,
{\it Geom. Funct. Anal.}, 14 (2004), 303-318.
\bibitem{hind} R. Hind, Lagrangian unknottedness in Stein surfaces, preprint math.SG/0311093.
\bibitem{hindivrii} R. Hind and A. Ivrii, Ruled 4-manifolds and
isotopies of symplectic surfaces, preprint math.SG/0510108.
\bibitem{hof} H. Hofer, K. Wysocki and E. Zehnder, A characterisation of the tight three-sphere, {\it Duke Math. J.}, 81(1995), no.1, 159-226.
\bibitem{hofa} H. Hofer, K. Wysocki and E. Zehnder, Properties of pseudoholomorphic curves in symplectisations I: Asymptotics, {\it Ann. Inst. H. Poincar\'{e} Anal. Non Lineaire}, 13(1996), no.3, 337-379.
\bibitem{hofi} H. Hofer, K. Wysocki and E. Zehnder, Properties of pseudoholomorphic curves in symplectisations II: Embedding controls and algebraic invariants, {\it Geom. Funct. Anal.}, 5(1995), no.2, 337-379.
\bibitem{hoff} H. Hofer, K. Wysocki and E. Zehnder, Properties of pseudoholomorphic curves in symplectisations III: Fredholm theory, {\it Topics in nonlinear analysis}, 381-475, Prog. Nonlinear Differential Equations Appl., 35, Birkh\"{a}user, Basel, 1999.
\bibitem{HLS}
H.~Hofer, V.~Lizan, and J.-C. Sikorav.
\newblock On genericity for holomorphic curves in four-dimensional
  almost-complex manifolds.
\newblock {\em J. Geom. Anal.}, 7(1):149--159, 1997.

\bibitem{ivrii} A. Ivrii.
\newblock Lagrangian unknottedness of tori in certain symplectic $4$-manifolds, Phd thesis.
\bibitem{wein} D. McDuff and D. Salamon, Introduction to symplectic topology. Second edition. Oxford Mathematical Monographs. The Clarendon Press, Oxford University Press, New York, 1998
\bibitem{MS2}
D.~McDuff and D.~Salamon.
\newblock {\em {$J$}-holomorphic curves and quantum cohomology}, volume~6 of
  {\em University Lecture Series}.
\newblock American Mathematical Society, Providence, RI, 1994.
\bibitem{msa} D. McDuff and D. Salamon, $J$-holomorphic curves and
symplectic topology. American Mathematical Society Colloquium
Publications, 52. American Mathematical Society, Providence, RI,
2004.

\bibitem{park} D. Park, M. Poddar, S. Vidussi, Homologous
Non-isotopic symplectic surfaces of higher genus, to appear in
Trans. Amer. Math. Soc.
\bibitem{rs} J. Robbin and D. Salamon, The Maslov index for paths, {\it Topology}, 32(1993), 827-844.

\bibitem{seidel} P. Seidel, Lagrangian two-spheres can be symplectically knotted, {\it J. Diff. Geom.}, 52(1999), 145-171.
\bibitem{ziller} W. Ziller, Geometry of the Katok examples, {\it Ergod. Th. and Dynam. Sys.}, 3(1982), 135-157.

\end{thebibliography}
\end{document}